\pdfoutput=1

\documentclass[a4paper,11pt]{article}
\usepackage[english]{babel}
\usepackage{amsfonts,amsmath,amssymb,theorem}
\usepackage{hyperref} 
\usepackage{colortbl}
\usepackage{ifthen}
\usepackage{hhline}

\makeatletter
\def\@seccntformat#1{\csname the#1\endcsname.\ } 
\def\@biblabel#1{#1.} 
\makeatother

\date{}

\newenvironment{proof}[1][\hspace{-1.0ex}]%
{\par\addvspace{1mm}{\sc Proof\hspace{1.0ex}{#1}.} }%
{\quad$\blacktriangle$\par\addvspace{1mm}}

\newif\ifNoRemark
\def\addtheorem#1#2#3#4{
\ifthenelse{\equal{#2}{}}{}%
{\ifthenelse{\expandafter\isundefined\csname the#2\endcsname}{\newcounter{#2}}{}}
\newenvironment{#1}[1][\global\NoRemarktrue]
{\par\addvspace{2mm plus 0.5mm minus 0.2mm}\noindent 
{\bf #3}\ifthenelse{\equal{#2}{}}{}%
{\refstepcounter{#2}{\bf ~\csname the#2\endcsname}}%
{\bf \vphantom{##1}\ifNoRemark.\ \else\ (##1).\fi}\begingroup #4}%
{\endgroup\par\addvspace{1mm plus 0.5mm minus 0.2mm}\global\NoRemarkfalse}
\expandafter\newcommand\csname b#1\endcsname{\begin{#1}}
\expandafter\newcommand\csname e#1\endcsname{\end{#1}}
}

\addtheorem{theorem}{thrm}{Theorem}{\sl}
\addtheorem{lemma}{lmm}{Lemma}{\sl}
\addtheorem{pro}{prpstn}{Proposition}{\sl}
\addtheorem{problem}{problem}{Problem}{\sl}

\begin{document}

\title{On the number of $n$-ary quasigroups of finite order\footnote{This work
was partially supported by the Federal Target Program
``Scientific and Educational Personnel of Innovation Russia'' for 2009-2013 (contract No.
02.740.11.0429) and the Russian Foundation for Basic Research
(grants 08-01-00671, 08-01-00673).}}
\author{V.\,N.\,Potapov$^{\star}$, D.\,S.\,Krotov$^{\star}$}

 \maketitle

\begin{center}
$^\star$\textit{Sobolev Institute of Mathematics, prospect Akademika Koptyuga 4, Novosibirsk\\
Mechanics and Mathematics Department of the
Novosibirsk State University, Pirogova 2, Novosibirsk }
\end{center}
\begin{abstract}

Let $Q(n,k)$ be the number of $n$-ary quasigroups of order $k$. We
derive a recurrent formula for $Q(n,4)$. We prove that for all
$n\geq 2$ and $k\geq 5$ the following inequalities hold:
$\left(\frac{k-3}2\right)^{\frac n2}\left(\frac{k-1}2\right)^{\frac
n2} < \log_2 Q(n,k) \leq {c_k(k-2)^{n}} $, where $c_k$ does not
depend on $n$. So, the upper asymptotic  bound for $Q(n,k)$ is
improved for any $k\geq 5$ and the lower bound is improved for odd
$k\geq 7$.

\end{abstract}

UDC: 519.143


\section{Introduction}

\hspace*{\parindent}
An algebraic system from a set $\Sigma$ of cardinality $|\Sigma|=k$ and $n$-ary operation $f:\Sigma^n\rightarrow \Sigma$ is called an
$n$-ary quasigroup of order $k$ iff the unary operation obtained by fixing any $n-1$ arguments of $f$ by any values from $\Sigma$ is always bijective.
The corresponding function $f$
is often also called an \emph{$n$-ary quasigroup} 
(the value table of such a function is known as {\em latin hypercube}; if $n=2$, {\em latin square}).

Let us fix the set $\Sigma=\{0,1,\dots,k-1\}$. Denote by $Q(n,k)$
the number of different $n$-ary quasigroups of order $k$ (for fixed
$\Sigma$)\footnote{Sometimes by the number of quasigroups one mean
the number of mutually nonisomorphic quasigroups.}. It is known that
for every $n$ there exist only two $n$-ary quasigroups of order $2$.
There are exactly $Q(n,3)=3\cdot2^n$ different $n$-ary quasigroups
of order $3$, which form one equivalence class. In
\cite{EN:PotKro:asymp} it is proved that $Q(n,4)=3^{n+1}2^{2^n
+1}(1+o(1))$ as $n\rightarrow\infty$. In Section~\ref{s:order4} we
suggest a recurrent way to calculate the numbers $Q(n,4)$ and list
the first $8$ values. Before, only five values of $Q(n,4)$ where
known; furthermore the numbers $Q(n,5)$ and $Q(n,6)$ are known for
$n\leq 5$ and $n\leq 3$ respectively, see~\cite{EN:MK-W:small},
and the number $Q(2,k)$ for $k\leq 11$, see \cite{EN:MK-W:11} and references there.

The asymptotics of the number and even of the logarithm of the number
(and even of the logarithm of the logarithm of the number) of
$n$-ary quasigroups of orders more than $4$ is unknown.
In \cite{EN:KPS:ir}, the following lower bounds are derived: $Q(n,5) \geq
2^{3^{n/3}-c}$, where $c<0.072$; $Q(n,k)\geq 2^{(k/2)^n}$ if $k$ is even;
$Q(n,k)\geq 2^{n(k/3)^n}$ if $k\equiv 0\bmod 3$; $Q(n,k)\geq
2^{1.5\lfloor k/3\rfloor^n}$ for an arbitrary $k$.
The following upper bound was established in \cite{EN:Pot:2008:LB}:
   $Q(n,k) \leq 3^{{(k-2)}^n}2^{n(k-2)^{n-1}} $.

In this paper we improve the upper bound (Section~\ref{s:UB}) on the number of $n$-ary
quasigroups of finite order and the lower bound (Section~\ref{s:LB}) on the number of $n$-ary
quasigroups of odd order:
$$
\left(\frac{k-3}2\right)^{\frac n2}\left(\frac{k-1}2\right)^{\frac n2} <
\log_2 Q(n,k) \leq {c_k(k-2)^{n}}, $$ where $c_k$ does not depend on $n$;
explicitly, $c_k=\frac{\log_2 k!}{k-2}+\frac{k}{k-4}$.

\section{An upper bound} \label{s:UB}

\hspace*{\parindent} We will say that a set $M\subseteq
\Sigma^n$ satisfies Property (A) iff for every element
$\bar x \in M$ and every position $i=1,\ldots,n$ there is 
another element $\bar y\in M$ differing from $\bar x$ only in the position $i$.
By induction it is easy to get the following:
\begin{pro}\label{chpro3}
Any nonempty subset $C\subseteq \Sigma^n$ that satisfies Property (A)
has the cardinality at least $2^n$.
\end{pro}

A function $g:\Omega\rightarrow\Sigma$ where $\Omega\subset
\Sigma^n$ is called a {\it partial $n$-ary quasigroup of order} $|\Sigma|$
if $g(\overline{x})\neq
g(\overline{y})$ for any two tuples  $\overline{x},
\overline{y}\in\Omega$ differing in exactly one position.
We will say that an $n$-ary quasigroup $f:\Sigma^n\rightarrow\Sigma$ is an \emph{extension}
of a partial $n$-ary quasigroup $g:\Omega\rightarrow\Sigma$ where $\Omega\subset \Sigma^n$
if $f|_{\Omega}\equiv g$.

 \begin{lemma}\label{chlem1} Let
$|\Sigma|=k$, $B=\Sigma\setminus\{a,b\}$, $k\geq 3$, $a,b\in
\Sigma$. Then a partial $n$-ary quasigroup $g:\Sigma^{n-1}\times B \rightarrow \Sigma$
has at most $2^{(k/2)^{n-1}}$ different extensions. \end{lemma}

\begin{proof}
Denote by $P$ the set of the unordered pairs of elements of $\Sigma$.
Consider a partial $n$-ary quasigroup $g:\Sigma^{n-1}\times B
\rightarrow \Sigma$. Define the function $G:\Sigma^{n-1}\rightarrow P$
by the equality $G(\overline{x})= \Sigma\setminus \{g(\overline{x}c) :
c\in \Sigma\setminus\{a,b\}\}$. Define the graph $\Gamma=\langle
\Sigma^{n-1}, E\rangle$ where two vertices $\overline{x}$ and
$\overline{y}$ are adjacent if and only if the tuples
$\overline{x}$ and $\overline{y}$ differ in exactly one position and
$G(\overline{x})\cap G(\overline{y})\neq \varnothing$.
It is easy to see that connected components of $\Gamma$
satisfy Property~(A).

Let $n$-ary quasigroups $f_1$ and $f_2$ be extensions of $g$. 
It is not difficult to see that $\{f_1(\overline{x}a), f_1(\overline{x}b)\}= G(\overline{x})$
for every $\overline{x}\in\Sigma^{n-1}$; moreover, if
$f_1(\overline{x}a)= f_2(\overline{x}a)$, then $f_1$ and
$f_2$ coincide on the whole connected component of $\Gamma$ containing $\overline{x}\in\Sigma^{n-1}$. 
So, to define an extension of $g$ uniquely it is sufficient to choose one from 
the two possible values for every connected component of $\Gamma$.
It follows from Proposition~\ref{chpro3} that every connected component has cardinality at least $2^{n-1}$. 
Then the number of connected components of $\Gamma$ does dot exceed ${(k/2)^{n-1}}$. 
Hence $g$ has not more than $2^{(k/2)^{n-1}}$ extensions.
\end{proof}

\begin{theorem} If $k\geq 5$ and $n\geq 2$ then $Q(n,k) \leq
2^{c_k(k-2)^{n}} $, where $c_k=\frac{\log_2 k!}{k-2}+\frac{k}{k-4}$.
\end{theorem}\begin{proof} The number of partial $n$-ary quasigroups
$g:\Sigma^{n-1}\times B \rightarrow \Sigma$, where $|\Sigma|=k$,
$B=\Sigma\setminus\{a,b\}$, does not exceed $Q(n,k)^{k-2}$. 
From Lemma~\ref{chlem1} we have
\begin{equation}\label{che1}
Q(n+1,k)\leq Q(n,k)^{k-2}2^{(k/2)^{n}}.
\end{equation}
Denote $\alpha_n=\log_2 Q(n,k)/(k-2)^n$. Then from
(\ref{che1}) we have  $$\alpha_{n+1}\leq \alpha_n
+\left(\frac{k}{2(k-2)}\right)^n.$$ Since $\alpha_1=\frac{\log_2
k!}{k-2}$ and
$\sum\limits_{n=1}^\infty\left(\frac{k}{2(k-2)}\right)^n=
\frac{k}{k-4}$, we get $\alpha_n\leq \frac{\log_2
k!}{k-2}+\frac{k}{k-4}$. \end{proof}

\section{A lower bound} \label{s:LB}

\hspace*{\parindent} Let $a$ and $b$ be two different elements of $\Sigma$. 
By \emph{$\{a,b\}$-com\-po\-nent} of an $n$-ary quasigroup$f$
we will mean such a set $S\subset\Sigma^n$ that 1)
$f(S)=\{a,b\}$ and 2) the function
$$
g(\bar x)=\begin{cases}
  f(\bar x) & \mbox{whenever }\bar x\not\in S, \cr
  b & \mbox{whenever }\bar x\in S \mbox{ and } f(\bar x)=a, \cr
  a & \mbox{whenever }\bar x\in S \mbox{ and } f(\bar x)=b
\end{cases}
$$
is also an $n$-ary quasigroup. In this case we will say that $g$ is obtained from $f$ by \emph{switching} the component $S$.
We note that in the definition of an $\{a,b\}$-com\-po\-nent the condition 2) can be replaced by Property~(A) from the previous section.
It is obvious that switching disjoint components can be performed independently:

\begin{pro}\label{chpro0}
Let $S$ and $S'$ be disjoint $\{a,b\}$- and $\{c,d\}$- (respectively) components of an $n$-ary quasigroup $f$.
Let an $n$-ary quasigroup $g$ is obtained from $f$ by switching $S$. Then $S'$ is a $\{c,d\}$-com\-po\-nent of $g$ too.
\end{pro}

The following proposition can be easily derived from the definition of an $\{a,b\}$-com\-po\-nent; similar statement can be found in
\cite{EN:KPS:ir}.

\begin{pro}\label{chpro5} Let $C=\{c_1,d_1\}\times
\{c_2,d_2\}$ be an $\{a,b\}$-com\-po\-nent of a $2$-ary qua\-si\-group $g$.
Let $C_i$ be a $\{c_i,d_i\}$-com\-po\-nent of an $n_i$-ary quasigroup $q_i$, $i=1,2$. 
Then the set $C_1\times C_2$ is an $\{a,b\}$-com\-po\-nent of the $(n_1+n_2)$-ary quasigroup
$f$, where $f(\bar x_1 ,\bar x_2 )\equiv g(q_1(\bar x_1),q_2( \bar x_2 ))$.\end{pro}

A $2$-ary quasigroup $\varphi:\Sigma\rightarrow\Sigma$ is called {\it idempotent} iff $\varphi(x,x)=x$ for every $x\in \Sigma$.
It is known (see, e.g., \cite{EN:Belousov:loops}) that
\begin{pro}\label{chpro4} For every $m\geq 3$ there exists an idempotent $2$-ary quasigroup of order $m$. \end{pro}

The following proposition presents a construction of $2$-ary quasigroups, which will be used
to establish a lower bound on the number of $n$-ary quasigroups of odd order.

\begin{pro}\label{chpro6} For any $m\geq 3$
there exists a $2$-ary quasigroup $\psi$ of order $2m+1$ that has
$m$ $\{2i,2i+1\}$-com\-po\-nents for every $i\in \{0,\dots,
m-1\}$; moreover, all except one $\{2i,2i+1\}$-com\-po\-nents
have form $\{2j,2j+1\}\times\{2l,2l+1\}$. \end{pro}

\begin{proof}
By Proposition~\ref{chpro4} there exists an idempotent
$2$-ary quasigroup $\varphi_m$ of order $m$. For each
$a,b\in\{0,\dots,m-1\}$, $a\neq b$, and $\delta,\sigma\in \{0,1\}$
define

$\psi(2a+\delta,2b+\sigma)= 2\varphi_m(a,b) + (\delta+\sigma\bmod 2);$

$\psi(2a+\delta,2a+\delta)= 2a+1-\delta;$

$\psi(2a+\delta,2a+1-\delta)= k-1;$

$\psi(k-1,2a+\delta)=\psi(2a+\delta,k-1)=2a+\delta;$

$\psi(k-1,k-1)=k-1.$

Straightforwardly, $\psi$ is a $2$-ary quasigroup that satisfied the desired 
properties.
 \end{proof}

The following is examples of the value tables of a $2$-ary quasigroup $\varphi_4$ 
and the corresponding $\psi$: 
\vskip5mm

$\varphi_4:$
\newcommand\numberWidth{0.5mm}
{
\newcommand\0[2]{\multicolumn{1}{#1>{\columncolor[rgb]{1.0,0.8,0.8}}c#2}{\makebox[\numberWidth][c]0}}
\newcommand\1[2]{\multicolumn{1}{#1>{\columncolor[rgb]{0.7,1.0,0.7}}c#2}{\makebox[\numberWidth][c]1}}
\newcommand\2[2]{\multicolumn{1}{#1>{\columncolor[rgb]{0.8,0.8,1.0}}c#2}{\makebox[\numberWidth][c]2}}
\newcommand\3[2]{\multicolumn{1}{#1>{\columncolor[rgb]{1.0,1.0,0.0}}c#2}{\makebox[\numberWidth][c]3}}
\begin{tabular}{cccc}
\hline
 \0{|}{}&\2{}{}&\3{}{}&\1{}{|}\\
 \3{|}{}&\1{}{}&\0{}{}&\2{}{|}\\
 \1{|}{}&\3{}{}&\2{}{}&\0{}{|}\\
 \2{|}{}&\0{}{}&\1{}{}&\3{}{|}\\
\hline
\end{tabular}
} \qquad $\psi:$
{
\newcommand\0[2]{\multicolumn{1}{#1>{\columncolor[rgb]{1.0,0.8,0.8}}c#2}{\makebox[\numberWidth][c]0}}
\newcommand\1[2]{\multicolumn{1}{#1>{\columncolor[rgb]{1.0,0.9,0.9}}c#2}{\makebox[\numberWidth][c]1}}
\newcommand\2[2]{\multicolumn{1}{#1>{\columncolor[rgb]{0.7,1.0,0.7}}c#2}{\makebox[\numberWidth][c]2}}
\newcommand\3[2]{\multicolumn{1}{#1>{\columncolor[rgb]{0.9,1.0,0.9}}c#2}{\makebox[\numberWidth][c]3}}
\newcommand\4[2]{\multicolumn{1}{#1>{\columncolor[rgb]{0.8,0.8,1.0}}c#2}{\makebox[\numberWidth][c]4}}
\newcommand\5[2]{\multicolumn{1}{#1>{\columncolor[rgb]{0.9,0.9,1.0}}c#2}{\makebox[\numberWidth][c]5}}
\newcommand\6[2]{\multicolumn{1}{#1>{\columncolor[rgb]{1.0,1.0,0.0}}c#2}{\makebox[\numberWidth][c]6}}
\newcommand\7[2]{\multicolumn{1}{#1>{\columncolor[rgb]{1.0,1.0,0.7}}c#2}{\makebox[\numberWidth][c]7}}
\newcommand\8[2]{\multicolumn{1}{#1>{\columncolor[rgb]{1.0,1.0,1.0}}c#2}{\makebox[\numberWidth][c]8}}
\setlength\arrayrulewidth{0.3pt}
\begin{tabular}{ccccccccc}
\hline
\1{|}{ }&\8{ }{ }&\4{|}{ }&\5{ }{ }&\6{|}{ }&\7{ }{ }&\2{|}{ }&\3{ }{|}&\0{ }{|}\\
\8{|}{ }&\0{ }{ }&\5{|}{ }&\4{ }{ }&\7{|}{ }&\6{ }{ }&\3{|}{ }&\2{ }{|}&\1{ }{|}\\
\hhline{|--------|~|}
\6{|}{ }&\7{ }{ }&\3{|}{ }&\8{ }{ }&\0{|}{ }&\1{ }{ }&\4{|}{ }&\5{ }{|}&\2{ }{|}\\
\7{|}{ }&\6{ }{ }&\8{|}{ }&\2{ }{ }&\1{|}{ }&\0{ }{ }&\5{|}{ }&\4{ }{|}&\3{ }{|}\\
\hhline{|--------|~|}
\2{|}{ }&\3{ }{ }&\6{|}{ }&\7{ }{ }&\5{|}{ }&\8{ }{ }&\0{|}{ }&\1{ }{|}&\4{ }{|}\\
\3{|}{ }&\2{ }{ }&\7{|}{ }&\6{ }{ }&\8{|}{ }&\4{ }{ }&\1{|}{ }&\0{ }{|}&\5{ }{|}\\
\hhline{|--------|~|}
\4{|}{ }&\5{ }{ }&\0{|}{ }&\1{ }{ }&\2{|}{ }&\3{ }{ }&\7{|}{ }&\8{ }{|}&\6{ }{|}\\
\5{|}{ }&\4{ }{ }&\1{|}{ }&\0{ }{ }&\3{|}{ }&\2{ }{ }&\8{|}{ }&\6{ }{|}&\7{ }{|}\\
\hhline{|--------|~|}
\0{|}{ }&\1{ }{ }&\2{ }{ }&\3{ }{ }&\4{ }{ }&\5{ }{ }&\6{ }{ }&\7{ }{ }&\8{ }{|}\\
\hline
\end{tabular}
}

From Proposition~\ref{chpro3} it is easy to conclude that the odd-order 
$2$-ary quasigroup constructed in Proposition~\ref{chpro6} has the maximal number of
mutually disjoint components among all $2$-ary quasigroups of the same order.

\begin{theorem} If $k$ is an odd integer ,$k\geq 5$, and $n\geq 2$, then
$$Q(n,k) \geq
2^{\left(\frac {k-3}{2}\right)^{\left\lfloor \frac{n-1}2\right\rfloor}\left(\frac {k-1}{2}\right)^{ \left\lceil\frac{n+1}2\right\rceil}} >
2^{\left(\frac {k-3}{2}\right)^{n/2}\left(\frac {k-1}{2}\right)^{n/2}} .$$
\end{theorem}

\begin{proof}
Let $\psi$ be the $2$-ary quasigroup of order $k$ constructed in Proposition~\ref{chpro6}.
 Define the $n$-ary quasigroup $\Psi^n$ by the following recurrent equalities:

$\Psi^2\equiv\psi$;

$\Psi^{2m+1}(\overline{x},y)=\psi(\Psi^{2m}(\overline{x}),y)$;

$\Psi^{2m+2}(\overline{x},y,z)=\psi(\Psi^{2m}(\overline{x}),\psi(y,z))$.

Denote by $\alpha_n$ the number of $\{2i,2i+1\}$-components
of $\Psi^{n}$ where $i\in\{0,\dots,\frac{k-3}{2}\}$. 
From Propositions~\ref{chpro5} and~\ref{chpro6} we have the relations 
$\alpha_{2}= \frac{k-1}{2}$,
$\alpha_{2m +1}\geq \alpha_{2m}\frac{k-3}{2}$, $\alpha_{2m +2}\geq
\alpha_{2m}\frac{k-3}{2}\frac{k-1}{2}$. Then $\alpha_{2m}\geq
\left(\frac{k-3}{2}\right)^{m-1}\left(\frac{k-1}{2}\right)^{m}$
and
$\alpha_{2m +1}\geq
\left(\frac{k-3}{2}\right)^{m}\left(\frac{k-1}{2}\right)^{m}$.

Since $\{2i,2i+1\}$-components with different $i$ are disjoint, 
the number of disjoint components is at least $\frac{k-1}2\alpha_n$. 
From Proposition~\ref{chpro0} we deduce that we can get the desired 
number of different $n$-ary quasigroups of order $k$ by switching disjoint 
components in $\Psi^n$.
\end{proof}

\section{The number of different $n$-ary quasigroups of order $4$} \label{s:order4}

\hspace*{\parindent}  Denote $[n]=\{1,\dots,n\}$.
An $n$-ary quasigroup $f$ is called an {\em $n$-ary loop} iff there exists
an element $e\in \Sigma$, which is called an {\em identity},
such that for all $i\in [n]$ and $a\in \Sigma$ it is true
$f(e\dots e\underset{i}{a}e\dots e) =a$. 
In what follows we always assume that $0$ is an identity 
of an $n$-ary loop
(in general, an $n$-ary loop can have more than one identities provided $n\geq3$). 
Especially we note that this agreement is essential in the treatment 
of the concept of the number of $n$-ary loops. 
In particular, we have the following simple and well-known fact:

\bpro \label{pro:12a} Let $Q'(n,k)$ be the number of $n$-ary loops of order $k$. Then \\ $Q(n,k)=k\cdot ((k-1)!)^n Q'(n,k)$. \epro

An $n$-ary quasigroop $f$ is called {\em permutably reducible} (we will omit ``permutably'') iff there exist an integer $m$, $2\leq m < n$,\quad an $(n-m+1)$-ary
quasigroup $h$,\quad an $m$-ary quasigroup $g$\quad, and a permutation
$\sigma: [n] \to [n]$  such that
$$f(x_1,\dots,x_{n}) \equiv h\left(g(x_{\sigma(1)},\dots, x_{\sigma(m)}),
 x_{\sigma(m+1)},\dots, x_{\sigma(n)}\right).$$

In this section we will assume that $\Sigma = \{0,1,2,3\}$;
i.e., we will consider only the $n$-ary quasigroups of order $4$.
It is known (see, e.g., \cite{EN:Belousov:loops}) that there are exactly four
binary loops of order $4$ (one is isomorphic to the group $Z_2\times Z_2$
and three, to the group $Z_4$).

The following statement is straightforward from the theorem in \cite{EN:Cher}.

\blemma \label{proCher} Every reducible $n$-ary loop $f$ of order $4$
 admits exactly one of the following two representations.
\begin{equation}
f(\overline{x})= q_0(q_1(\tilde x_1),...,q_m(\tilde
x_m))\label{eq:Decomp.-of-quas.}
\end{equation}
 where 
 $q_j$ are $n_j$-ary loops;
 $ \tilde x_j$ are tuples of variables $x_i$, $i\in I_j$, where
 $\{I_j\}$ is a partition of $[n]$, $j=1,\ldots, m$; 
 $q_0$ is an irreducible $m$-ary loop, $m\geq 3$. 
 Moreover, the partition $\{I_j\}$ in such a representation is unique for every $f$.

 \begin{equation}
f(\overline{x})=q_1(\tilde x_1)\ast...\ast q_k(\tilde
x_k) \label{eq:Decomp.-of-quas1.}
\end{equation}
where $\ast$ is a binary operation in one of the $4$ loops, $q_j$, $j=1,\ldots,k$, are
$n_j$-ary loops that are not representable in the form
$q_j(\tilde x_j) = q'(\tilde x'_j)\ast q''(\tilde x''_j)$,
 $ \tilde x_j$ are tuples of variables $x_i$, $i\in I_j$, where
 $\{I_j\}$ is a partition of $[n]$. 
 Moreover, the partition $\{I_j\}$ in such a representation is unique for every $f$.
\elemma

By the {\em root operation} of an $n$-ary quasigroup $f$ 
we will mean the $m$-ary quasigroup $q_0$ if 
 (\ref{eq:Decomp.-of-quas.}) holds, and the binary operation $\ast$
if (\ref{eq:Decomp.-of-quas1.}) holds.

Simple combinatorial calculations give the following formula for the number
$F_{\overline{j},\overline{k}}$ of different partitions of $[n]$
into $k$ subsets from which exactly $k_i$ subsets have cardinality
$j_i$, $1\leq i \leq t$, $0<j_1<\dots<j_t$:

\begin{equation}
F_{\overline{j},\overline{k}}=\frac{n!}{(j_1!)^{k_1}\dots(j_t!)^{k_t}}\frac{1}{{k_1!}\dots
{k_t!}},\label{coeff}
\end{equation}
where $k_1+k_2+\dots+k_t=k$, $k_1j_1+k_2j_2+\dots+k_tj_t=n$.

Let $f:\Sigma^n\rightarrow \Sigma$ be an $n$-ary quasigroup;
define the set 
$$ S_{a,b}(f) \triangleq \cup\{\bar x\in
\Sigma^n: f(\bar x)\in\{a,b\}\}.$$ 
An $n$-ary loop $f$ will be called {\em $a$-semilinear}, where
$a\in\{1,2,3\}$, if the characteristic function $\chi_{S_{0,a}(f)}$ of the set $S=S_{0,a}(f)$ has the form
\begin{equation}\label{eq:def_lin}
\chi_{S_{0,a}(f)}(x_1,\dots,x_{n}) \equiv \sum_{i=1}^n \chi_{\{0,a\}}(x_i) \bmod 2.
\end{equation}
An $n$-ary loop $f$ is called {\em linear} if it is $a$-semi\-linear and $b$-semi\-linear for some different $a$ and $b$ from
$\{1,2,3\}$. It is not difficult to check the following:
\bpro \label{pro14} One of the four binary loops of order $4$ is linear
(the one that is isomorphic to $Z_2\times Z_2$); the other three are 
$1$-, $2$-, and $3$- semilinear respectively.\epro

It is known (see \cite{EN:PotKro:asymp}) that
\bpro \label{pro:Linear-Quas.-Is-Alone} A linear $n$-ary loop is unique and is
$1$-, $2$-, and $3$- semilinear.
\epro

It is not difficult to see (see also \cite{EN:PotKro:asymp}) the following:
 \bpro \label{pro:13b} Let $f$ be a reducible
 $a$-semi\-linear $n$-ary loop; then $f$ can be represented as the composition
{\rm(\ref{eq:Decomp.-of-quas.})} or
{\rm(\ref{eq:Decomp.-of-quas1.})} of $a$-semi\-linear loops.
 \epro

 Let us denote by $\ell^a_n$ the number of the
$a$-semi\-linear $n$-ary loops and by $\ell_n$ 
the number of the semilinear $n$-ary loops.

As established in \cite{EN:PotKro:asymp},
the number of the $n$-ary loops asymptotically coincides with $\ell_n$, which can be easily calculated:
\blemma[\cite{EN:PotKro:asymp}]
\label{th:Numb.-of-Semilin.-Quas.} $\ell_n=3\cdot 2^{2^n-n-1} -2$,
$\ell^a_n= 2^{2^n-n-1}$ for $a\in \{1,2,3\}$. \elemma

In \cite{EN:KroPot:4} the set of $n$-ary quasigroups of order $4$  was characterized 
in the terms defined above; namely, the following was proved:

\btheorem\label{top} Every $n$-ary loop of order $4$ is reducible or semilinear. \etheorem

This fact gives a base for deriving a recurrent formula for the number of $n$-ary loops (and quasigroups) of order $4$.

We will use the following notation:

$v_n$ is the number of the $n$-ary loops (of order $4$);

$r_n^*$ is the number of the $n$-ary loops with the binary root operation $*$ (all of them are reducible except $*$ itself);

$r_n^0$ is the number of the reducible $n$-ary loops with the root operation of arity at least $3$;

$r_n^{a*}$ is the number of the $a$-semi\-linear $n$-ary loops with the
$a$-semi\-linear binary root operation $*$ (all of them are reducible except $*$ itself);

$r_n^{a0}$ is the number of the reducible $a$-semi\-linear $n$-ary loops with the
root operation of arity at least $3$;

$p_n^a$ is the number of irreducible $a$-semi\-linear $n$-ary loops;

$p_n$ is the number of irreducible $n$-ary loops.

From Lemma~\ref{proCher} and Proposition~\ref{pro:13b}, the following relations follow:
\def\sumKJ{\sum_{\bar j, \bar k}}
$$r_n^{a*}= \sum\limits_{i=2}^n \sumKJ F_{\overline{j},\overline{k}}(\ell^a_{j_1}
-r_{j_1}^{a*})^{k_1}\cdots(\ell^a_{j_t} -r_{j_t}^{a*})^{k_t},$$
$$r_n^{*}= \sum\limits_{i=2}^n \sumKJ F_{\overline{j},\overline{k}}(v_{j_1}
-r_{j_1}^{*})^{k_1}\cdots(v_{j_t} -r_{j_t}^{*})^{k_t},$$
$$r_n^{a0}= \sum\limits_{i=3}^{n-1}p_i^a \sumKJ F_{\overline{j},\overline{k}}(\ell^a_{j_1})^{k_1}\cdots(\ell^a_{j_t})^{k_t},$$
$$r_n^{0}= \sum\limits_{i=3}^{n-1}p_i \sumKJ F_{\overline{j},\overline{k}}(v_{j_1})^{k_1}\cdots(v_{j_t})^{k_t},$$
where the second sum os over the tuples $\bar k =(k_1,\ldots,k_t)$ and
$\bar j =(j_1,\ldots,j_t)$ of positive integers satisfying $k_1+\dots+k_t=i$,
$k_1j_1+k_2j_2+\dots+k_tj_t=n$ and $j_1<\dots<j_t$. From
Theorem~\ref{top} and Proposition~\ref{pro:Linear-Quas.-Is-Alone}
we have $v_n= p_n+r_n^0 + 4r^*_n$,
$p^a_n=\ell^a_n-r_n^{a0}-2r_n^{a*}$, $p_n=3p^a_n$. From Lemma~\ref{th:Numb.-of-Semilin.-Quas.},  $\ell^a_n= 2^{2^n-n-1}$ for
$a\in \{1,2,3\}$.

The initial values are trivial: $r_1^{a*}=r_1^{*}=r_1^{a0}=r_1^0=0$. 
We see that the equalities above and Proposition~\ref{pro:12a} provide a recurrent way of calculation
of the number of the $n$-ary quasigroups of order $4$.

Finally, we list the first eight values of $Q'(n,4)$:\\
1,\\
4,\\
64,\\
7132,\\
201538000,\\
432345572694417712,\\
3987683987354747642922773353963277968,\\
{678469272874899582559986240285280710364867063489779510427038722229750276832},\\
and of $Q(n,4)$:\\
24,\\
576,\\
55296,\\
36972288,\\
6268637952000,\\
80686060158523011084288,\\
4465185218736554544676917926460256725000192,\\
{4558271384916189349044295395852008182480786230841798008741684281906576963885826048}.\\
The program on \texttt{python} to calculate these numbers can be found in the appendix.

\section{Conclusion}
\hspace*{\parindent} 
We will briefly discuss 
a connection of our topic with the known concept of latin traid.
A partial $n$-ary quasigroup $t: \Omega \to
\Sigma$, $\Omega\subset\Sigma^n$ is called
a \emph{multidimensional latin trade}, here for briefity simply \emph{trade}
iff there exist another partial $n$-ary quasigroup $t': \Omega
\to \Sigma$ such that

1) $t(\bar x) \neq t'(\bar x)$ for all $\bar x\in \Omega$;

2) for any $i$ from $1$ to $n$ and any admissible values
$x_1$, \ldots, $x_{i-1}$, $x_{i-1}$, \ldots, $x_n$ the sets
$\{t(x_1, \ldots, x_{i-1}, y, x_{i-1}, \ldots, x_n) \ |\ y\in\Sigma
\}$ and $\{t'(x_1, \ldots, x_{i-1}, y, x_{i-1}, \ldots, x_n) \ |\
y\in\Sigma \}$ coincide.

In this case the pair $(t,t')$ is called a \emph{bitrade} 
(depending on the context, bitrades are considered either as ordered,
or as unordered pairs); the trade $t'$ is called a {\em mate} of $t$. 
In the case $n=2$ bitrades
(latin bitrades) are widely studied, see the survey
\cite{EN:Cav:rev}.

We will say that an $n$-ary quasigroup $f$ has a trade $t$
iff $t=f|_{\Omega}$ for some $\Omega$. 
As follows from the definitions,
replacing of the values of $f$ in $\Omega$
by the values of a mate $t'$ of $t$ results in another $n$-ary
quasigroup. 
We will say that trades $t=f|_{\Omega}$ and $s=f|_{\Theta}$ are 
\emph{independent} iff their supports $\Omega$
and $\Theta$ are disjoint. The maximal number of mutually independent trades
of an $n$-ary quasigroup $f$ will be called its \emph{trade number} $\mathrm{trd}(f)$. 
Denote by $\mathrm{Trd}(n,k)$ the maximum of ${\rm trd}(f)$ over the all $n$-ary quasigroups $f$ of order $k$. Since independent trades of an $n$-ary quasigroup
can be independently replaced by mates, the number $Q(n,k)$
of different $n$-ary quasigroups of order $k$ satisfies
\begin{equation}\label{eq:Trn}
 Q(n,k)\geq 2^{{\rm Trd}(n,k)}.
\end{equation}
It is easy to understand that the lower bound in Section~\ref{s:LB} 
(as well as all the bounds in \cite{EN:KPS:ir}) is derived by this way:
$\{a,b\}$-com\-po\-nent is the support of some trade by definitions. 
Since the support of a trade satisfies Property~(A), Proposition~\ref{chpro3} 
implies ${\rm Trd}(n,k)\leq
k^n/2^n=2^{(\log_2 k-1)n}$; moreover, for even $k$ the equality is easily proved.
For odd $k$, as follows from the results of Section~\ref{s:LB}, we have
${\rm Trd}(n,k)\geq 2^{c(k)n}$ where $c(k)
\underset{k\to\infty}{\longrightarrow} \log_2k-1$. 
But for fixed $k$, in particular, for the small values $5$, $7$, \ldots,
the question about the asymptotics of ${\rm Trd}(n,k)$
remains open.

\begin{problem}
Establish the asymptotics of the logarithm and the asymptotics of the value ${\rm Trd}(n,k)$ as $n\to\infty$ for odd $k\geq 5$.
\end{problem}

Another question is how the estimation (\ref{eq:Trn}) close to the real value. 
For the order $4$ it is asymptotically tight in logarithms.  
For any larger fixed order the asymptotics of $\log\log Q(n,k)$ is unknown.
It is natural to conjecture that the asymptotics of $\log\log Q(n,k)$ and $\log{\rm Trd}(n,k)$ coincide.

\begin{problem}
  Is it true that ${\displaystyle\lim_{n\to\infty}}\left(\frac{\log_2\log_2Q(n,k)}n\right)={\displaystyle\lim_{n\to\infty}}\left(\frac{\log_2{\rm Trd}(n,k)}n\right)$?
  In particular, is it true that ${\displaystyle\lim_{n\to\infty}}\left(\frac{\log_2\log_2Q(n,k)}n\right) \leq \log_2 k -1$?
\nopagebreak\end{problem}
Even the existence of these limits is not proved yet.


\providecommand\href[2]{#2} \providecommand\url[1]{\href{#1}{#1}} \def\DOI#1{
  {DOI}: \href{http://dx.doi.org/#1}{#1}}\def\DOIURL#1#2{ {DOI}:
  \href{http://dx.doi.org/#2}{#1}}\providecommand\bbljun{June}

\newpage
\appendix
\section{Appendix. The numbers $Q(n,4)$ and $Q'(n,4)$: program}
\begin{verbatim}
# PYTHON program
# do not forget leading spaces!
N=12 \# the maximum arity to calculate
P=[[{}]]+[[] for x in [0]*N]
for n in range(N):
  for Po in P[n]:
    for j in range(N-n-sum(Po.values())+1):
      Pn={i+1:Po[i] for i in Po}
      if j>0: Pn.update({1:j})
      if n+j>0: P[n+sum(Po.values())+j]+=[Pn]

\# now P[n] is a list of partitions of n into positive summands:
\# E.g., 16==1*3+2*4+5*1==(1+1+1)+(2+2+2+2)+(5) is saved as P[16][51]=={1:3, 2:4, 5:1}

for Pn in P: Pn.pop()  \# remove trivial 1-partitions

from math import factorial
mul = lambda a,b:a*b
Aut = lambda Pt:reduce(mul,[factorial(Pt[nu])*factorial(nu)**Pt[nu] for nu in Pt], 1)
F=[[factorial(n)/Aut(Pt) for Pt in P[n]] for n in range(N+1)]
\# F[n][i] is the number of partitions of an n-set 
\# that correspond to the partition P[n][i] of n.
La=[2L**(2**n-n-1) for n in range(N+1)]
Ras,Ra0,R_0,R_s,P_a,V,T = [0,0L],[0, 0L],[0,0L],[0,0L],[0,0L],[1,1L],[4,24L]
for n in range(2, N+1):
  V+=[0L]; T+=[0L]; P_a+=[0L]; Ras+=[0L]; Ra0+=[0L]; R_0+=[0L]; R_s+=[0L]
  for i in range(len(P[n])):
    R_0[n], Ra0[n], R_s[n], Ras[n] = map(lambda A, B, C :
      A[n] + reduce(mul,[(B[nu]-C*A[nu])**P[n][i][nu] for nu in P[n][i]],1) 
       *((1-C)*P_a[sum(P[n][i].values())]+C)*F[n][i],
        (R_0, Ra0, R_s, Ras), (V, La, V, La), (0, 0, 1, 1))
  R_0[n] *= 3
  P_a[n] = La[n] - Ra0[n] - 2*Ras[n]
  V[n] = 3*P_a[n] + R_0[n] + 4*R_s[n]
  T[n] = 4*(6**n)*V[n]

print "\n Reduced (A211214):", V
print "\n Total (A211215):", T 
\end{verbatim}

\end{document}


\title{О числе {$n$}"=арных квазигрупп конечного порядка
\footnote{Работа
выполнена при поддержке ФЦП <<Научные и научно-педагогические кадры
инновационной России>> на 2009--2013 гг. (гос. контракт №
02.740.11.0429) и Российского фонда фундаментальных исследований
(проекты 08-01-00671, 08-01-00673).}}
\author{В.\,Н.\,Потапов$^{\star}$, Д.\,С.\,Кротов$^{\star}$}

 \maketitle

\begin{center}
$^\star$\textit{Институт математики им. С.\,Л.\,Соболева СО РАН,\\ проспект Академика Коптюга 4, Новосибирск\\
Механико-математический факультет, Новосибирский государственный университет, ул. Пирогова 2, Новосибирск }
\end{center}
\begin{abstract}
Пусть $Q(n,k)$ --- число $n$-арных квазигрупп порядка $k$. Получена
рекуррентная формула для чисел $Q(n,4)$. Доказано, что при любых
$n\geq 2$ и $k\geq 5$ справедливы неравенства $\left(\frac{k-3}2\right)^{\frac n2}\left(\frac{k-1}2\right)^{\frac n2} < \log_2 Q(n,k) \leq
{c_k(k-2)^{n}} $, где $c_k$ не зависит от $n$. Таким образом,
верхняя асимтотическая граница для чисел $Q(n,k)$ улучшена при любых
$k\geq 5$, нижняя "--- при нечётных $k\geq 7$.

{\bf Ключевые слова}: $n$-арная квазигруппа, латинский куб, лупа,
асимптотика, компонента, латинский трэйд.
\end{abstract}

УДК: 519.143

MSC2010: 20N15, 05B15

\section{Введение}

\hspace*{\parindent}  Алгебраическая система из множества $\Sigma$
мощности $|\Sigma|=k$ и $n$-арной операции \linebreak $f:
\Sigma^n\rightarrow \Sigma$ называется $n$-арной квазигруппой
порядка $k$, если унарная операция, полученная фиксацией любых $n-1$
аргументов операции $f$ любыми значениями из $\Sigma$, всегда
является биекцией. Принято называть {\it $n$-арной квазигруппой
порядка} $k$  или {\it $n$-квазигруппой порядка} $k$ также и
соответствующую функцию $f$ (таблица значений такой функции известна
под названием {\em латинский гиперкуб}, в случае $n=2$ "--- {\em
латинский квадрат}).

Зафиксируем множество элементов $\Sigma=\{0,1,\dots,k-1\}$.
Обозначим через $Q(n,k)$ число различных $n$-арных квазигрупп порядка $k$
(при фиксированном $\Sigma$)\footnote{Иногда под числом квазигрупп
 подразумевают число попарно неизоморфных квазигрупп.}.
Известно, что для  каждого $n$ существует только две $n$-арных
квазигруппы порядка $2$. Имеется $Q(n,3)=3\cdot2^n$ различных
$n$-арных квазигрупп порядка $3$, которые составляют единственный
класс эквивалентности. В \cite{PotKro:asymp.ru} доказано асимптотическое
равенство $Q(n,4)=3^{n+1}2^{2^n +1}(1+o(1))$ при
$n\rightarrow\infty$.  В настоящей работе (раздел~\ref{s:order4}) предложен рекуррентный
способ вычисления чисел $Q(n,4)$, выписаны первые $8$ значений этой величины.
Ранее были известны только пять первых значений $Q(n,4)$,
кроме того, известны числа $Q(n,5)$ и $Q(n,6)$ для $n\leq 5$ и $n\leq 3$ соответственно, см.~\cite{MK-W:small},
и число $Q(2,k)$ при $k\leq 11$, см. работу \cite{MK-W:11} и
библиографию в ней.

 Асимптотика числа и даже логарифма числа
(и даже логарифма логарифма числа) $n$"=арных квазигрупп порядков
больше $4$ неизвестна. В \cite{KPS:ir} получены следующие нижние
оценки $Q(n,5) \geq 2^{3^{n/3}-c}$, где $c<0.072$; $Q(n,k)\geq
2^{(k/2)^n}$  если $k$ чётно;  $Q(n,k)\geq 2^{n(k/3)^n}$ если $k$
кратно трём; $Q(n,k)\geq 2^{1.5\lfloor k/3\rfloor^n}$ для
произвольного $k$.
  В
\cite{Pot:2008:LB.ru} была предложена  верхняя оценка числа $n$"=арных
квазигрупп
 порядка $k$:
   $Q(n,k) \leq 3^{{(k-2)}^n}2^{n(k-2)^{n-1}} $.

В настоящей статье усилена верхняя оценка (раздел~\ref{s:UB}) числа $n$"=арных
квазигрупп
 конечного порядка  и нижняя оценка (раздел~\ref{s:LB}) числа
$n$"=арных квазигрупп
 нечётного порядка:
$$
\left(\frac{k-3}2\right)^{\frac n2}\left(\frac{k-1}2\right)^{\frac n2} <
\log_2 Q(n,k) \leq {c_k(k-2)^{n}}, $$ где $c_k$ не зависит от $n$, точнее, $c_k=\frac{\log_2 k!}{k-2}+\frac{k}{k-4}$.

\section{Верхняя оценка} \label{s:UB}

\hspace*{\parindent} Будем говорить, что множество $M\subseteq
\Sigma^n$ удовлетворяет свойству (A), если для любого элемента
$\bar x$ из $M$ и каждой позиции $i=1,\ldots,n$ найдется другой элемент $\bar y$ из $M$, отличающийся от $\bar x$ только в позиции $i$. По индукции легко получить
следующее:
\begin{pro}\label{chpro3}
Любое непустое подмножество $C\subseteq \Sigma^n$,
 удовлетворяющее свойству (A),
 имеет мощность не менее $2^n$.
\end{pro}

{\it Частичной $n$"=арной квазигруппой порядка} $|\Sigma|$
называется функция $g:\Omega\rightarrow\Sigma$, где $\Omega\subset
\Sigma^n$, удовлетворяющая следующему свойству: $g(\overline{x})\neq
g(\overline{y})$ для  любых двух наборов $\overline{x},
\overline{y}\in\Omega$, отличающихся ровно в одной позиции. Будем
говорить, что
$n$"=арная квазигруппа $f:\Sigma^n\rightarrow\Sigma$ является \emph{продолжением}
частичной $n$"=арной квазигруппы $g:\Omega\rightarrow\Sigma$, где $\Omega\subset \Sigma^n$,
если $f|_{\Omega}\equiv g$.

 \begin{lemma}\label{chlem1} Пусть
$|\Sigma|=k$, $B=\Sigma\setminus\{a,b\}$, $k\geq 3$, $a,b\in
\Sigma$. Тогда частичная $n$"=арная квазигруппа
$g:\Sigma^{n-1}\times B \rightarrow \Sigma$ имеет не более чем
$2^{(k/2)^{n-1}}$ различных продолжений. \end{lemma}

\begin{proof}

 Пусть $P$ ---
множество неупорядоченных пар элементов множества $\Sigma$.
Рассмотрим частичную $n$"=арную квазигруппу $g:\Sigma^{n-1}\times B
\rightarrow \Sigma$. Определим функцию $G:\Sigma^{n-1}\rightarrow P$
равенством $G(\overline{x})= \Sigma\setminus \{g(\overline{x}c) :
c\in \Sigma\setminus\{a,b\}\}$. Определим граф $\Gamma=\langle
\Sigma^{n-1}, E\rangle$, в котором две вершины $\overline{x}$ и
$\overline{y}$ соединены ребром тогда и только тогда, когда наборы
$\overline{x}$ и $\overline{y}$ отличаются только в одной позиции и
$G(\overline{x})\cap G(\overline{y})\neq \varnothing$. Нетрудно
видеть, что компоненты связности графа $\Gamma$ удовлетворяют
свойству (A).

 Пусть $n$"=арные квазигруппы $f_1$ и $f_2$ являются
продолжениями частичной $n$"=арной квазигруппы $g$. Нетрудно видеть,
что $\{f_1(\overline{x}a), f_1(\overline{x}b)\}= G(\overline{x})$
для любого $\overline{x}\in\Sigma^{n-1}$, причём если
$f_1(\overline{x}a)= f_2(\overline{x}a)$, то продолжения $f_1$ и
$f_2$ совпадают на всей компоненте связности графа $\Gamma$,
содержащей вершину $\overline{x}\in\Sigma^{n-1}$. Таким образом, для
однозначного определения продолжения частичной $n$"=арной
квазигруппы $g$ достаточно  зафиксировать одно из двух возможных
значений в каждой компоненте связности графа $\Gamma$. Из
предложения \ref{chpro3} следует, что каждая  компонента связности
имеет мощность не  менее $2^{n-1}$. Тогда число компонент связности
графа $\Gamma$ не превосходит ${(k/2)^{n-1}}$. А значит, $g$ имеет не более
$2^{(k/2)^{n-1}}$ продолжений.
\end{proof}

\begin{theorem} Если $k\geq 5$ и $n\geq 2$, то $Q(n,k) \leq
2^{c_k(k-2)^{n}} $, где $c_k=\frac{\log_2 k!}{k-2}+\frac{k}{k-4}$.
\end{theorem}\begin{proof} Число частичных $n$"=арных квазигрупп
$g:\Sigma^{n-1}\times B \rightarrow \Sigma$, где $|\Sigma|=k$,
$B=\Sigma\setminus\{a,b\}$ не превосходит $Q(n,k)^{k-2}$. Из леммы
\ref{chlem1} следует неравенство
\begin{equation}\label{che1}
Q(n+1,k)\leq Q(n,k)^{k-2}2^{(k/2)^{n}}.
\end{equation}
Введём обозначение $\alpha_n=\log_2 Q(n,k)/(k-2)^n$. Тогда из
неравенства (\ref{che1}) имеем  $$\alpha_{n+1}\leq \alpha_n
+\left(\frac{k}{2(k-2)}\right)^n.$$ Поскольку $\alpha_1=\frac{\log_2
k!}{k-2}$ и
$\sum\limits_{n=1}^\infty\left(\frac{k}{2(k-2)}\right)^n=
\frac{k}{k-4}$, имеем $\alpha_n\leq \frac{\log_2
k!}{k-2}+\frac{k}{k-4}$. \end{proof}

\section{Нижняя оценка} \label{s:LB}

\hspace*{\parindent} Пусть $a$ и $b$ --- два различных элемента из
$\Sigma$. \emph{$\{a,b\}$"=Компонентой} $n$"=арной квазигруппы $f$
будем называть такое множество $S\subset\Sigma^n$, что 1)
$f(S)=\{a,b\}$ и 2) функция
$$
g(\bar x)=\begin{cases}
  f(\bar x) & \mbox{при }\bar x\not\in S, \cr
  b & \mbox{при }\bar x\in S \mbox{ и } f(\bar x)=a, \cr
  a & \mbox{при }\bar x\in S \mbox{ и } f(\bar x)=b.
\end{cases}
$$
 также является $n$"=арной квазигруппой. В этом случае будем говорить, что $g$ получена из $f$ \emph{свитчингом} компоненты $S$.
Заметим, что в определении $\{a,b\}$"=компоненты условие 2) можно заменить свойством (A) из предыдущего раздела.
Очевидно, что свитчинг непересекающихся компонент можно производить независимо:

\begin{pro}\label{chpro0}
Пусть $S$ и $S'$ "--- непересекающиеся $\{a,b\}$- и $\{c,d\}$- (соответственно) компоненты $n$-арной квазигруппы $f$ и $n$-арная квазигруппа $g$ получена из $f$ свитчингом компоненты $S$. Тогда $S'$ также является $\{c,d\}$-компонентой квазигруппы $g$.
\end{pro}

Следующее предложение нетрудно получить из
определения $\{a,b\}$-компоненты, аналогичное утверждение имеется в
\cite{KPS:ir}.

\begin{pro}\label{chpro5} Пусть множество $C=\{c_1,d_1\}\times
\{c_2,d_2\}$ является $\{a,b\}$-компонентой $2$-квазигруппы $g$.
Пусть множество $C_i$ является $\{c_i,d_i\}$-компонентой
$n_i$"=арной квазигруппы $q_i$ при $i=1,2$. Тогда множество
$C_1\times C_2$ является $\{a,b\}$-компонентой $(n_1+n_2)$"=арной
квазигруппы $f$, где $f(\bar x_1 ,\bar x_2 )\equiv g(q_1(\bar x_1
),q_2( \bar x_2 ))$.\end{pro}

$2$-Квазигруппа $\varphi:\Sigma\rightarrow\Sigma$  называется {\it
идемпотентной}, если $\varphi(x,x)=x$ для любого $x\in \Sigma$.
Известно (см., например, \cite{Belousov:loops.ru}), что верно
\begin{pro}\label{chpro4} Для любого $m\geq 3$ имеется
идемпотентная $2$"=квазигруппа
 порядка $m$. \end{pro}

В следующем предложении приведена  конструкция $2$"=квазигрупп,
которая будет использована при доказательстве нижней оценки числа
$n$"=арных квазигрупп нечётного порядка.

\begin{pro}\label{chpro6} Для любого  $m\geq 3$
найдётся $2$"=квазигруппа $\psi$ порядка $2m+1$,  имеющая
$m$ $\{2i,2i+1\}$-компонент для каждого $i\in \{0,\dots,
m-1\}$, причём все кроме одной $\{2i,2i+1\}$-компоненты
имеют вид $\{2j,2j+1\}\times\{2l,2l+1\}$. \end{pro}

\begin{proof}
По предложению~\ref{chpro4} найдётся идемпотентная
$2$"=квазигруппа $\varphi_m$ порядка $m$. Для любых
$a,b\in\{0,\dots,m-1\}$, $a\neq b$, и $\delta,\sigma\in \{0,1\}$
определим

$\psi(2a+\delta,2b+\sigma)= 2\varphi_m(a,b) + (\delta+\sigma\bmod 2);$

$\psi(2a+\delta,2a+\delta)= 2a+1-\delta;$

$\psi(2a+\delta,2a+1-\delta)= k-1;$

$\psi(k-1,2a+\delta)=\psi(2a+\delta,k-1)=2a+\delta;$

$\psi(k-1,k-1)=k-1.$

Непосредственная проверка показывает, что $\psi$ есть $2$"=квазигруппа,
обладающая требуемыми свойствами.
 \end{proof}

Ниже приведён пример таблиц значений $2$"=квазигруппы $\varphi_4$ и соответствующей $\psi$: \vskip5mm

$\varphi_4:$
\newcommand\numberWidth{0.5mm}
{
\newcommand\0[2]{\multicolumn{1}{#1>{\columncolor[rgb]{1.0,0.8,0.8}}c#2}{\makebox[\numberWidth][c]0}}
\newcommand\1[2]{\multicolumn{1}{#1>{\columncolor[rgb]{0.7,1.0,0.7}}c#2}{\makebox[\numberWidth][c]1}}
\newcommand\2[2]{\multicolumn{1}{#1>{\columncolor[rgb]{0.8,0.8,1.0}}c#2}{\makebox[\numberWidth][c]2}}
\newcommand\3[2]{\multicolumn{1}{#1>{\columncolor[rgb]{1.0,1.0,0.0}}c#2}{\makebox[\numberWidth][c]3}}
\begin{tabular}{cccc}
\hline
 \0{|}{}&\2{}{}&\3{}{}&\1{}{|}\\
 \3{|}{}&\1{}{}&\0{}{}&\2{}{|}\\
 \1{|}{}&\3{}{}&\2{}{}&\0{}{|}\\
 \2{|}{}&\0{}{}&\1{}{}&\3{}{|}\\
\hline
\end{tabular}
} \qquad $\psi:$
{
\newcommand\0[2]{\multicolumn{1}{#1>{\columncolor[rgb]{1.0,0.8,0.8}}c#2}{\makebox[\numberWidth][c]0}}
\newcommand\1[2]{\multicolumn{1}{#1>{\columncolor[rgb]{1.0,0.9,0.9}}c#2}{\makebox[\numberWidth][c]1}}
\newcommand\2[2]{\multicolumn{1}{#1>{\columncolor[rgb]{0.7,1.0,0.7}}c#2}{\makebox[\numberWidth][c]2}}
\newcommand\3[2]{\multicolumn{1}{#1>{\columncolor[rgb]{0.9,1.0,0.9}}c#2}{\makebox[\numberWidth][c]3}}
\newcommand\4[2]{\multicolumn{1}{#1>{\columncolor[rgb]{0.8,0.8,1.0}}c#2}{\makebox[\numberWidth][c]4}}
\newcommand\5[2]{\multicolumn{1}{#1>{\columncolor[rgb]{0.9,0.9,1.0}}c#2}{\makebox[\numberWidth][c]5}}
\newcommand\6[2]{\multicolumn{1}{#1>{\columncolor[rgb]{1.0,1.0,0.0}}c#2}{\makebox[\numberWidth][c]6}}
\newcommand\7[2]{\multicolumn{1}{#1>{\columncolor[rgb]{1.0,1.0,0.7}}c#2}{\makebox[\numberWidth][c]7}}
\newcommand\8[2]{\multicolumn{1}{#1>{\columncolor[rgb]{1.0,1.0,1.0}}c#2}{\makebox[\numberWidth][c]8}}
\setlength\arrayrulewidth{0.3pt}
\begin{tabular}{ccccccccc}
\hline
\1{|}{ }&\8{ }{ }&\4{|}{ }&\5{ }{ }&\6{|}{ }&\7{ }{ }&\2{|}{ }&\3{ }{|}&\0{ }{|}\\
\8{|}{ }&\0{ }{ }&\5{|}{ }&\4{ }{ }&\7{|}{ }&\6{ }{ }&\3{|}{ }&\2{ }{|}&\1{ }{|}\\
\hhline{|--------|~|}
\6{|}{ }&\7{ }{ }&\3{|}{ }&\8{ }{ }&\0{|}{ }&\1{ }{ }&\4{|}{ }&\5{ }{|}&\2{ }{|}\\
\7{|}{ }&\6{ }{ }&\8{|}{ }&\2{ }{ }&\1{|}{ }&\0{ }{ }&\5{|}{ }&\4{ }{|}&\3{ }{|}\\
\hhline{|--------|~|}
\2{|}{ }&\3{ }{ }&\6{|}{ }&\7{ }{ }&\5{|}{ }&\8{ }{ }&\0{|}{ }&\1{ }{|}&\4{ }{|}\\
\3{|}{ }&\2{ }{ }&\7{|}{ }&\6{ }{ }&\8{|}{ }&\4{ }{ }&\1{|}{ }&\0{ }{|}&\5{ }{|}\\
\hhline{|--------|~|}
\4{|}{ }&\5{ }{ }&\0{|}{ }&\1{ }{ }&\2{|}{ }&\3{ }{ }&\7{|}{ }&\8{ }{|}&\6{ }{|}\\
\5{|}{ }&\4{ }{ }&\1{|}{ }&\0{ }{ }&\3{|}{ }&\2{ }{ }&\8{|}{ }&\6{ }{|}&\7{ }{|}\\
\hhline{|--------|~|}
\0{|}{ }&\1{ }{ }&\2{ }{ }&\3{ }{ }&\4{ }{ }&\5{ }{ }&\6{ }{ }&\7{ }{ }&\8{ }{|}\\
\hline
\end{tabular}
}

Из предложения \ref{chpro3} нетрудно заключить, что $2$"=квазигруппа
нечётного порядка $k$, построенная в предложении \ref{chpro6}, имеет
максимальное число непересекающихся компонент среди всех
$2$"=квазигрупп порядка $k$.

\begin{theorem} Если $k\geq 5$ --- нечётное и $n\geq 2$, то
$$Q(n,k) \geq
2^{\left(\frac {k-3}{2}\right)^{\left\lfloor \frac{n-1}2\right\rfloor}\left(\frac {k-1}{2}\right)^{ \left\lceil\frac{n+1}2\right\rceil}} >
2^{\left(\frac {k-3}{2}\right)^{n/2}\left(\frac {k-1}{2}\right)^{n/2}} .$$
\end{theorem}

\begin{proof}
Пусть $\psi$ -- $2$-квазигруппа порядка $k$, построенная в
предложении~\ref{chpro6}.
 Определим рекуррентно $n$"=арную квазигруппу $\Psi^n$
равенствами:

$\Psi^2\equiv\psi$;

$\Psi^{2m+1}(\overline{x},y)=\psi(\Psi^{2m}(\overline{x}),y)$;

$\Psi^{2m+2}(\overline{x},y,z)=\psi(\Psi^{2m}(\overline{x}),\psi(y,z))$.

Обозначим через $\alpha_n$ --- число $\{2i,2i+1\}$-компонент
$n$"=арной квазигруппы $\Psi^{n}$, где $i\in
\{0,\dots,\frac{k-3}{2}\}$. Из предложений \ref{chpro5} и
\ref{chpro6} имеем соотношения $\alpha_{2}= \frac{k-1}{2}$,
$\alpha_{2m +1}\geq \alpha_{2m}\frac{k-3}{2}$, $\alpha_{2m +2}\geq
\alpha_{2m}\frac{k-3}{2}\frac{k-1}{2}$. Тогда $\alpha_{2m}\geq
\left(\frac{k-3}{2}\right)^{m-1}\left(\frac{k-1}{2}\right)^{m}$
и
$\alpha_{2m +1}\geq
\left(\frac{k-3}{2}\right)^{m}\left(\frac{k-1}{2}\right)^{m}$.

Поскольку $\{2i,2i+1\}$-компоненты при различных $i$ не
пересекаются, всего непересекающихся компонент не меньше, чем
$\frac{k-1}2\alpha_n$. Из предложения \ref{chpro0} следует, что  из
$n$"=арной квазигруппы $\Psi^n$ свитчингами непересекающихся
компонент можно получить требуемое число различных $n$"=арных
квазигрупп порядка $k$.
\end{proof}

\section{Число $n$-арных квазигрупп порядка $4$} \label{s:order4}

\hspace*{\parindent}  Введём обозначение $[n]=\{1,\dots,n\}$.
$n$-Арная квазигруппа $f$ называется  {\em $n$-арной лупой}, если
существует такой элемент $e\in \Sigma$, называемый {\em единичным},
что для всех $i\in [n]$ и $a\in \Sigma$ имеет место равенство
$f(e\dots e\underset{i}{a}e\dots e) =a$. Далее мы всегда будем
подразумевать, что $0$ является единичным элементом $n$-арной лупы
(в общем случае могут быть и другие единичные элементы). Особо
отметим, что данное соглашение существенно в интерпретации понятия
числа $n$-арных луп. В частности, имеем следующий простой и хорошо
известный факт:

\bpro \label{pro:12a} Пусть $Q'(n,k)$ "--- число  $n$"=арных луп
порядка $k$. Тогда \\ $Q(n,k)=k\cdot ((k-1)!)^n Q'(n,k)$. \epro

$n$-Арная квазигруппа $f$ называется {\em разделимой (приводимой)},
если найдутся: целое число $m$, $2\leq m < n$,\quad $(n-m+1)$-арная
квазигруппа $h$,\quad $m$"=арная квазигруппа $g$\quad и перестановка
$\sigma: [n] \to [n]$  "--- такие, что
$$f(x_1,\dots,x_{n}) \equiv h\left(g(x_{\sigma(1)},\dots, x_{\sigma(m)}),
 x_{\sigma(m+1)},\dots, x_{\sigma(n)}\right).$$

В дальнейшем будем предполагать, что $\Sigma =
\{0,1,2,3\}$, т.\,е. речь пойдёт только об $n$"=арных квазигруппах
порядка $4$. Известно (см., например, \cite {Belousov:loops.ru}), что
имеется
ровно четыре бинарные лупы порядка $4$ (одна изоморфна группе $Z_2\times Z_2$
и три "--- группе $Z_4$).

 Следующее утверждение
является прямым следствием теоремы из \cite{Cher.ru}.

\blemma \label{proCher}  Для  разделимой $n$-арной лупы $f$ порядка
$4$
 справедливо одно и только одно из двух представлений.
\begin{equation}
f(\overline{x})= q_0(q_1(\tilde x_1),...,q_m(\tilde
x_m)),\label{eq:Decomp.-of-quas.}
\end{equation}
 где $q_j$ есть $n_j$-арные лупы при $j,1\leq j\leq m$, $q_0$ есть
 неразделимая $m$-арная лупа, $m\geq 3$,
 $ \tilde x_j$ "--- некоторые наборы переменных $x_i$, $i\in I_j$, где
 $\{I_j\}$ "--- разбиение
 множества $[n]$. Причём в данном представлении
 разбиение $\{I_j\}$ единственно.

 \begin{equation}
f(\overline{x})=q_1(\tilde x_1)\ast...\ast q_k(\tilde
x_k),\label{eq:Decomp.-of-quas1.}
\end{equation}
где $\ast$ есть бинарная операция в одной из $4$ луп, $q_j$ есть
$n_j$-арные лупы при $j, 1\leq j\leq k$, непредставимые в виде
$q_j(\tilde x_j) = q'(\tilde x'_j)\ast q''(\tilde x''_j)$,
 $ \tilde x_j$ "--- некоторые наборы переменных $x_i$, $i\in I_j$, где $\{I_j\}$ "--- разбиение
 множества $[n]$. Причём в данном представлении
 разбиение $\{I_j\}$ единственно.
\elemma

{\em Корневой операцией} $n$-арной квазигруппы $f$ будем называть
$m$-арную квазигруппу $q_0$, если имеет место представление
 (\ref{eq:Decomp.-of-quas.}), и бинарную операцию
$\ast$, если имеет место представление (\ref{eq:Decomp.-of-quas1.}).

Простой комбинаторный подсчёт показывает, что число
$F_{\overline{j},\overline{k}}$ различных разбиений множества $[n]$
на $k$ подмножеств, из которых $k_i$ подмножеств имеет мощность
$j_i$, $1\leq i \leq t$, $0<j_1<\dots<j_t$, удовлетворяет равенству

\begin{equation}
F_{\overline{j},\overline{k}}=\frac{n!}{(j_1!)^{k_1}\dots(j_t!)^{k_t}}\frac{1}{{k_1!}\dots
{k_t!}},\label{coeff}
\end{equation}
где $k_1+k_2+\dots+k_t=k$, $k_1j_1+k_2j_2+\dots+k_tj_t=n$.

Пусть $f:\Sigma^n\rightarrow \Sigma$ "--- $n$"=арная квазигруппа,
определим множество $$ S_{a,b}(f) \triangleq \{\bar x\in
\Sigma^n: f(\bar x)\in\{a,b\}\}.$$ $n$-Арную
лупу $f$ назовём {\em $a$"=полулинейной}, где
$a\in\{1,2,3\}$ если характеристическая функция $\chi_{S_{0,a}(f)}$ множества $S=S_{0,a}(f)$ имеет вид
\begin{equation}\label{eq:def_lin}
\chi_{S_{0,a}(f)}(x_1,\dots,x_{n}) \equiv \sum_{i=1}^n \chi_{\{0,a\}}(x_i) \bmod 2.
\end{equation}
$n$-Арная
лупа $f$ называется {\em линейной}, если она одновременно
является $a$"=полулинейной и $b$"=полулинейной, где $a\neq b$,
$a,b\in\{1,2,3\}$. Непосредственой проверкой нетрудно убедиться, что
справедливо
\bpro \label{pro14} Из четырёх бинарных луп порядка $4$ одна
(изоморфная группе $Z_2\times Z_2$) является линейной, а три
остальных $1$-, $2$- и $3$- полулинейными соответственно.\epro

 Известно (см. \cite{PotKro:asymp.ru}), что
\bpro \label{pro:Linear-Quas.-Is-Alone} Линейная $n$-арная
лупа единственна и является одновременно $1$-, $2$- и $3$- полулинейной.
\epro

Нетрудно видеть (см., также, \cite{PotKro:asymp.ru}), что справедливо
 \bpro \label{pro:13b} Пусть $f$ "--- разделимая
 $a$"=полулинейная
$n$-арная лупа, тогда $f$ можно представить как суперпозицию
$a$"=полулинейных луп вида {\rm(\ref{eq:Decomp.-of-quas.})} или
{\rm(\ref{eq:Decomp.-of-quas1.})}.\epro

 Обозначим через $\ell^a_n$ мощность множества
$a$"=полулинейных $n$-арных луп и через $\ell_n$ мощность множества
полулинейных $n$-арных луп.

Как было установлено в \cite{PotKro:asymp.ru},
мощность множества всех $n$-арных луп асимптотически
 совпадает с мощностью множества полулинейных $n$-арных луп, которая легко вычисляется:
\blemma[\cite{PotKro:asymp.ru}]
\label{th:Numb.-of-Semilin.-Quas.} $\ell_n=3\cdot 2^{2^n-n-1} -2$,
$\ell^a_n= 2^{2^n-n-1}$ при $a\in \{1,2,3\}$. \elemma

В   \cite {KroPot:4} получено описание $n$-арных квазигрупп порядка $4$ в
определённых выше терминах, а именно, доказана

\btheorem\label{top} Каждая $n$-арная  лупа  порядка $4$
является разделимой или полулинейной. \etheorem

На этом описании по-существу основывается вывод рекуррентной формулы
для числа $n$-арных луп (и квазигрупп) порядка $4$.

Введём следующие обозначения:

$v_n$ "--- число   $n$-арных луп (порядка $4$);

$r_n^*$ "--- число $n$-арных луп с бинарной корневой 
операцией $*$ (все такие лупы разделимы, 
за исключением собственно бинарной операции $*$);

$r_n^0$ "--- число  разделимых $n$-арных луп с корневой операцией
арности большей либо равной $3$;

$r_n^{a*}$ "--- число $a$"=полулинейных $n$-арных луп с
$a$"=полулинейной бинарной корневой операцией $*$ (все такие лупы разделимы, 
за исключением собственно бинарной операции $*$);

$r_n^{a0}$ "--- число  разделимых $a$"=полулинейных $n$-арных луп с
корневой операцией арности большей либо равной $3$;

$p_n^a$ "--- число  неразделимых $a$"=полулинейных $n$-арных луп;

$p_n$ "--- число  неразделимых  $n$-арных луп.

Из леммы \ref{proCher} и предложения \ref{pro:13b} вытекают
следующие соотношения:
\def\sumKJ{\sum_{\bar j, \bar k}}
$$r_n^{a*}= \sum\limits_{i=2}^n \sumKJ F_{\overline{j},\overline{k}}(\ell^a_{j_1}
-r_{j_1}^{a*})^{k_1}\cdots(\ell^a_{j_t} -r_{j_t}^{a*})^{k_t},$$
$$r_n^{*}= \sum\limits_{i=2}^n \sumKJ F_{\overline{j},\overline{k}}(v_{j_1}
-r_{j_1}^{*})^{k_1}\cdots(v_{j_t} -r_{j_t}^{*})^{k_t},$$
$$r_n^{a0}= \sum\limits_{i=3}^{n-1}p_i^a \sumKJ F_{\overline{j},\overline{k}}(\ell^a_{j_1})^{k_1}\cdots(\ell^a_{j_t})^{k_t},$$
$$r_n^{0}= \sum\limits_{i=3}^{n-1}p_i \sumKJ F_{\overline{j},\overline{k}}(v_{j_1})^{k_1}\cdots(v_{j_t})^{k_t},$$
где вторая сумма берётся по наборам $\bar k =(k_1,\ldots,k_t)$ и
$\bar j =(j_1,\ldots,j_t)$ положительных целых чисел,
удовлетворяющим равенствам $k_1+\dots+k_t=i$,
$k_1j_1+k_2j_2+\dots+k_tj_t=n$ и неравенствам $j_1<\dots<j_t$. Из
теоремы \ref{top} и предложения \ref{pro:Linear-Quas.-Is-Alone}
вытекают соотношения $v_n= p_n+r_n^0 + 4r^*_n$,
$p^a_n=\ell^a_n-r_n^{a0}-2r_n^{a*}$, $p_n=3p^a_n$. Из леммы
\ref{th:Numb.-of-Semilin.-Quas.} имеем  $\ell^a_n= 2^{2^n-n-1}$ при
$a\in \{1,2,3\}$.

Начальные значения для
перечисленных выше величин тривиальны: $r_1^{a*}=r_1^{*}=r_1^{a0}=r_1^0=0$. 
Нетрудно видеть, что приведённые выше равенства
и предложение \ref{pro:12a} обеспечивают рекуррентный способ
вычисления числа $n$-арных квазигрупп порядка $4$.

Наконец, выпишем первые восемь значений величины $Q'(n,4)$:
1,
4,
64,\linebreak
7132,
201538000,
432345572694417712,
3987683987354747642922773353963277968,\linebreak
{678469272874899582559986240285280710364867063489779510427038722229750276832},~---\
и\linebreak величины $Q(n,4)$:
24,
576,
55296,
36972288,
6268637952000,\linebreak
80686060158523011084288,
4465185218736554544676917926460256725000192,\linebreak
{4558271384916189349044295395852008182480786230841798008741684281906576963885826048}.\\
Программа на \texttt{python} для вычисления этих чисел приведена в приложении.

\section{Заключение}
\hspace*{\parindent} В заключении скажем несколько слов о связи
тематики настоящей статьи с известным понятием латинского трэйда
(latin trade). Частичная $n$-арная квазигруппа $t: \Omega \to
\Sigma$,\linebreak $\Omega\subset\Sigma^n$ называется
\emph{многомерным латинским трэйдом}, далее просто \emph{трэйдом},
если существует другая частичная $n$-арная квазигруппа $t': \Omega
\to \Sigma$ такая, что

1) $t(\bar x) \neq t'(\bar x)$ для всех $\bar x\in \Omega$;

2) для любого $i$, $i=1,\dots,n$, и для любых  допустимых значений
$x_1$, \ldots, $x_{i-1}$, $x_{i-1}$, \ldots, $x_n$ множества
$\{t(x_1, \ldots, x_{i-1}, y, x_{i-1}, \ldots, x_n) \ |\ y\in\Sigma
\}$ и $\{t'(x_1, \ldots, x_{i-1}, y, x_{i-1}, \ldots, x_n) \ |\
y\in\Sigma \}$ совпадают.

В этом случае пара $(t,t')$ называется \emph{битрейдом} (в
зависимости от контекста, битрейд рассматривают  как упорядоченную,
либо как неупорядоченную пару), а трэйд $t'$ называется
\emph{партнёром} трэйда $t$. Исследованию битрэйдов в случае $n=2$
(латинских битрэйдов) уделяется значительное внимание, см. обзор
\cite{Cav:rev}.

Будем говорить, что $n$-арная квазигруппа $f$ содержит трэйд $t$,
если $t=f|_{\Omega}$ для некоторого $\Omega$. При этом из
определений следует, что замена значений $f$ на множестве $\Omega$
значениями партнёра $t'$ трэйда $t$ приводит к другой $n$-арной
квазигруппе. Трэйды $t=f|_{\Omega}$ и $s=f|_{\Theta}$ назовём
\emph{независимыми}, если их носители (области определения) $\Omega$
и $\Theta$ не пересекаются. Максимальное число попарно независимых
трэйдов, которые содержит $n$-арная квазигруппа $f$, назовём её
\emph{трэйдовым числом} ${\rm trd}(f)$. Максимум ${\rm trd}(f)$ по
всем $n$-арным квазигруппам $f$ порядка $k$ обозначим через ${\rm
Trd}(n,k)$. Поскольку независимые трэйды в $n$-арной квазигруппе
можно независимо заменять на партнёров, число $Q(n,k)$ различных
$n$-арных квазигрупп порядка $k$ удовлетворяет неравенству
\begin{equation}\label{eq:Trn}
 Q(n,k)\geq 2^{{\rm Trd}(n,k)}.
\end{equation}
Легко понять, что нижняя оценка в разделе~\ref{s:LB} (как и все
оценки в \cite{KPS:ir}) получена именно таким образом:
$\{a,b\}$"=компонента по определению является носителем некоторого
трэйда. Поскольку носитель трэйда обладает свойством (A), из
предложения~\ref{chpro3} вытекает ${\rm Trd}(n,k)\leq
k^n/2^n=2^{(\log_2 k-1)n}$, причем для чётных $k$ легко доказать
равенство. Для нечётных $k$ из результатов раздела~\ref{s:LB}
следует оценка ${\rm Trd}(n,k)\geq 2^{c(k)n}$, где $c(k)
\underset{k\to\infty}{\longrightarrow} \log_2k-1$. Однако для
фиксированных $k$, в частности, для малых значений $5$, $7$, \ldots
вопрос об асимптотическом поведении величины ${\rm Trd}(n,k)$
остаётся открытым.

\begin{problem}
  Вычислить асимптотику логарифма и асимптотику величины ${\rm Trd}(n,k)$ при $n\to\infty$ для нечётных $k\geq 5$.
\end{problem}

Другой вопрос состоит в том, насколько оценка  (\ref{eq:Trn}) близка
к истинной. Для порядка $4$ оценка (\ref{eq:Trn}) асимптотически
точна после логарифмирования.  Для большего фиксированного порядка
асимптотика двукратного  логарифма величины $Q(n,k)$ неизвестна.
Кажется  естественным предположить, что  асимптотика двукратного
логарифма  величины  $Q(n,k)$ и логарифма величины ${\rm Trd}(n,k)$
совпадают.

\begin{problem}
  Верно ли, что ${\displaystyle\lim_{n\to\infty}}\left(\frac{\log_2\log_2Q(n,k)}n\right)={\displaystyle\lim_{n\to\infty}}\left(\frac{\log_2{\rm Trd}(n,k)}n\right)$?
  В частности, верно ли, что ${\displaystyle\lim_{n\to\infty}}\left(\frac{\log_2\log_2Q(n,k)}n\right) \leq \log_2 k -1$?
\end{problem}
Существование этих пределов также не доказано.


\newpage
\appendix
\section{Приложение. Программа для вычисления $Q(n,4)$ и $Q'(n,4)$}
\begin{verbatim}
# PYTHON program
# do not forget leading spaces!
N=12 \# the maximum arity to calculate
P=[[{}]]+[[] for x in [0]*N]
for n in range(N):
  for Po in P[n]:
    for j in range(N-n-sum(Po.values())+1):
      Pn={i+1:Po[i] for i in Po}
      if j>0: Pn.update({1:j})
      if n+j>0: P[n+sum(Po.values())+j]+=[Pn]

\# now P[n] is a list of partitions of n into positive summands:
\# E.g., 16==1*3+2*4+5*1==(1+1+1)+(2+2+2+2)+(5) is saved as P[16][51]=={1:3, 2:4, 5:1}

for Pn in P: Pn.pop()  \# remove trivial 1-partitions

from math import factorial
mul = lambda a,b:a*b
Aut = lambda Pt:reduce(mul,[factorial(Pt[nu])*factorial(nu)**Pt[nu] for nu in Pt], 1)
F=[[factorial(n)/Aut(Pt) for Pt in P[n]] for n in range(N+1)]
\# F[n][i] is the number of partitions of an n-set 
\# that correspond to the partition P[n][i] of n.
La=[2L**(2**n-n-1) for n in range(N+1)]
Ras,Ra0,R_0,R_s,P_a,V,T = [0,0L],[0, 0L],[0,0L],[0,0L],[0,0L],[1,1L],[4,24L]
for n in range(2, N+1):
  V+=[0L]; T+=[0L]; P_a+=[0L]; Ras+=[0L]; Ra0+=[0L]; R_0+=[0L]; R_s+=[0L]
  for i in range(len(P[n])):
    R_0[n], Ra0[n], R_s[n], Ras[n] = map(lambda A, B, C :
      A[n] + reduce(mul,[(B[nu]-C*A[nu])**P[n][i][nu] for nu in P[n][i]],1) 
       *((1-C)*P_a[sum(P[n][i].values())]+C)*F[n][i],
        (R_0, Ra0, R_s, Ras), (V, La, V, La), (0, 0, 1, 1))
  R_0[n] *= 3
  P_a[n] = La[n] - Ra0[n] - 2*Ras[n]
  V[n] = 3*P_a[n] + R_0[n] + 4*R_s[n]
  T[n] = 4*(6**n)*V[n]

print "\n Reduced (A211214):", V
print "\n Total (A211215):", T 
\end{verbatim}